\documentclass[a4paper]{article}%
\usepackage{amsmath}
\usepackage{amsfonts}
\usepackage{amssymb}
\usepackage{graphicx}%
\providecommand{\U}[1]{\protect\rule{.1in}{.1in}}
\newtheorem{theorem}{Theorem}

\newtheorem{proposition}[theorem]{Proposition}
\newtheorem{remark}[theorem]{Remark}

\newenvironment{proof}[1][Proof]{\noindent\textbf{#1.} }{\ \rule{0.5em}{0.5em}}
\begin{document}

\title{Relation between Exponential Moment problem and Classical Moment problem}
\author{O. Kounchev, H. Render, Ts. Tsachev}
\maketitle

\begin{abstract}
The main result of the paper is an interesting relation between the solution
of the truncated Exponential Moment problem and truncated Classical Moment
problem, considered on the half-line or on a compact interval.

\textbf{AMS Classification, MSC2010}: 44A60, 30E05, 42A82.


\end{abstract}

\section{Introduction}

The main result of the present paper is a curious relation between the
Exponential Moment problem and the Classical Moment problem \cite{Akhiezer}.
The generalized Moment problem which contains the exponential moment problem
has been considered for the first time by A. Markov in his original paper of
$1886,$ in the form
\[%
{\displaystyle\int_{a}^{b}}
u_{j}\left(  t\right)  d\mu\left(  t\right)  =c_{j}%
\]
where $\left\{  u_{j}\right\}  _{j=1}^{n}$ is a set of functions which satisfy
some approprate conditions, see \cite{Markov-1886}, also results and
historical remarks in \cite{KreinNudelman}. For the particular case of
functions of the form
\[
u_{j}\left(  t\right)  =e^{\lambda_{j}t}=x^{\lambda_{j}}%
\]
(where $x=e^{t}$ ) this problem has been later studied by F. Hausdorff
(1921-1923) and S. Bernstein (1928), cf. results and historical remarks in
\cite{Akhiezer} and \cite{KreinNudelman}.

The interest to the exponential polynomials, i.e. linear combinations of
exponential functions $\left\{  e^{\lambda_{j}t}\right\}  $ for real and
complex numbers $\lambda_{j}$ is motivated by the fact that the polyharmonic
functions in special domains are splitted into sums which contain special
instances of exponential functions, cf. \cite{kounchevBOOK}. Recently, they
have found numerous applications in Approximation theory, Spline theory, CAGD,
Image processing, Moment problems, Data Analysis, and others, as may be seen
in the following references,
\cite{AldazKounchevRender-ConstructiveApproximation-2008}, \cite{NATO-Paper},
\cite{KoReArkiv2010}, \cite{KoRe13}, \cite{KoReInter2019},
\cite{KounchevRender2021JCAM}, \cite{KounchevRenderTsachevBIT}, \cite{Unser},
\cite{Gustafsson-Putinar}, \cite{McCa91}, \cite{Schumaker},
\cite{RamsaySilverman}.

The paper is organized as follows: in Section 2 we recall basic facts and
background material about exponential polynomials. In Section 3 we present and
prove our main results, Theorem \ref{THM-Main} and Theorem \ref{THM-Moments}.

\section{Exponential polynomials}

Exponential polynomials in a general setting have been given a special
attention e.g. in \cite{Berenstein}. In the present paper we will consider
exponential polynomials for some special parameters $\Lambda$.

For simplicity sake, in the present paper we will assume that the vector of
the "frequencies" $\Lambda=\left(  \lambda_{0},\lambda_{1},...,\lambda
_{N}\right)  $ (in the engineering literature frequencies are called only the
imaginary parts of the constants $\lambda_{j}$) consists of pair-wise
different real numbers $\lambda_{i}>0$ and $N$ is an \textbf{odd integer}. We
assume that $\mathbb{S}$ is a linear functional on the linear span
\[
U_{N}=\left\langle e^{\lambda_{0}x},....,e^{\lambda_{N}x}\right\rangle
_{\operatorname{lin}}\text{ with }N=2n-1,
\]
hence, $U_{N}$ contains an even number of elements. First of all, we consider
a special basis of the space $U_{N}$ defined by the set
\[
B\left(  x\right)  :=\left\{  b_{j}\left(  x\right)  \right\}  _{j=0}^{N}%
\]
where we have put
\begin{equation}
b_{j}\left(  x\right)  :=b_{N,j}\left(  x\right)  :=j!\cdot\Phi_{\left(
\lambda_{0},....,\lambda_{N}\right)  }^{\left(  N-j\right)  }\left(  x\right)
\text{ for }j=0,....,N.\label{bj-DEFINITION}%
\end{equation}

The function (sometimes called fundamental) $\Phi\left(  t\right)
:=\Phi_{\left(  \lambda_{0},....,\lambda_{N}\right)  }\left(  x\right)  $ is
defined as the unique element of the space $U_{N}$ which satisfies
\[
\Phi^{\left(  j\right)  }\left(  0\right)  =0\qquad\text{for }%
j=0,1,...,N-1;\qquad\Phi^{\left(  N\right)  }\left(  0\right)  =1.
\]

For the proof that the system of functions $B$ is a basis for the space
$U_{N}$ we refer to \cite{AldazKounchevRender-ConstructiveApproximation-2008},
\cite{NATO-Paper}. The "moment" linear functional $\mathbb{S}:U_{N}%
\rightarrow\mathbb{R}$ is defined by the finite number of "exponential
moments":
\begin{equation}
\widehat{c}_{j}:=\widehat{c}_{N,j}:=\int j!\cdot\Phi_{\left(  \lambda
_{0},....,\lambda_{N}\right)  }^{\left(  N-j\right)  }\left(  x\right)
d\mu\left(  x\right)  =\int b_{j}\left(  x\right)  d\mu\left(  x\right)
\text{ for }j=0,...,N. \label{cjHAT}%
\end{equation}

\section{Main results}

We prove the following main result:

\begin{theorem}
\label{THM-Main} For the function $\Phi_{N}\left(  x\right)  =\Phi_{\left(
\lambda_{0},...,\lambda_{N}\right)  }\left(  x\right)  ,$ for $x\in\mathbb{R}$
with $N$ odd, we assume that
\[
\Phi_{N}^{\left(  j\right)  }\left(  x\right)  \geq0\text{ for all }x>0\text{
and all }j\geq0.
\]
Then, the following two statements hold:

1. For every $x\in\mathbb{R}_{+}$ the sequence $\left\{  b_{j}\left(
x\right)  \right\}  _{j=0}^{N}$ defined in (\ref{bj-DEFINITION}) is
\textbf{non-negative} w.r.t. the half-axis $\mathbb{R}_{+}$, i.e. for every
$k$ with $2k\leq N,$ and for every non-zero vector $p\in\mathbb{R}^{k+1},$ the
inequality
\begin{equation}
Q_{1}\left(  p\right)  :=%
{\displaystyle\sum_{i,j=0}^{k}}
p_{i}p_{j}b_{i+j}\left(  x\right)  \geq0\label{POSITIVE-whole-AXIS}%
\end{equation}
holds, as well as the inequality
\begin{equation}
Q_{2}\left(  p\right)  :=%
{\displaystyle\sum_{i,j=0}^{k}}
p_{i}p_{j}b_{i+j+1}\left(  x\right)  \geq0\label{POSITIVE-HALF-AXIS}%
\end{equation}

2. For every $x\in\left[  0,1\right]  ,$ the sequence $\left\{  b_{j}\left(
x\right)  \right\}  _{j=0}^{N}$ is \textbf{positive } with respect to the
interval $\left[  0,1\right]  ,$ i.e. for every non-zero vector $p\in
\mathbb{R}^{k+1}$ with $2k+1\leq N,$ the two inequalities hold, namely,
inequality (\ref{POSITIVE-HALF-AXIS}), and second, inequality
\begin{equation}%
{\displaystyle\sum_{i,j=0}^{k}}
p_{i}p_{j}\left(  b_{i+j}\left(  x\right)  -b_{i+j+1}\left(  x\right)
\right)  \geq0. \label{POSITIVE-INTERVAL}%
\end{equation}

\end{theorem}

\begin{remark}
\label{REMARK-1} In the Theorem abve we use the definition of
\textbf{non-negativity} of a \textbf{sequence} with respect to the whole axis,
the half-axis and to a compact interval $\left[  a,b\right]  .$ We refer to
the classical monogrph of Akhiezer \cite{Akhiezer}, Definition $2.6.3,$ on p.
$70$ and next pages; The equivalent criteria which we use here below, are in
formulas $(2.50)$ in \cite{Akhiezer}, p. $74,$ and section $6.5$ on p. $76.$
See also Krein-Nudelman \cite{KreinNudelman}. Thus, the first statement of the
Theorem is in fact the criterion for solvability of the Stieltjes moment
problem (on the half-line), \cite{Akhiezer}, p. $76,$ formulas in section
$6.5.$ The second statement is the criterion for solvability of the Hausdorff
moment problem on a compact interval, see \cite{Akhiezer}, p. $74,$ formulas
$(2.50).$
\end{remark}

\begin{proof}
of Theorem \ref{THM-Main}:

As is proved in \cite{KoRe13}, also in \cite{NATO-Paper}, we can expand the
function $\Phi_{N}$ in a Taylor series, namely:
\[
\Phi_{N}\left(  x\right)  =%
{\displaystyle\sum_{s=N}^{\infty}}
a_{s}x^{s}\text{ with }a_{N}=\frac{1}{N!}.
\]
By our assumption we have $a_{s}\geq0$ for all $s=0,....$. Note that for every
$m$ with $0\leq m\leq N$ holds
\[
\Phi_{N}^{\left(  N-m\right)  }\left(  x\right)  =\frac{d^{N-m}}{dx^{N-m}}%
\Phi_{N}\left(  x\right)  =%
{\displaystyle\sum_{s=N}^{\infty}}
a_{s}\frac{d^{N-m}}{dx^{N-m}}x^{s},
\]
hence, the vector of functions $\left\{  b_{j}\left(  x\right)  \right\}
_{j=0}^{2k}$ is equal to
\[
A_{N,2k}\left(  x\right)  :=\left(
{\displaystyle\sum_{s=N}^{\infty}}
a_{s}\frac{d^{N}}{dx^{N}}x^{s},1!%
{\displaystyle\sum_{s=N}^{\infty}}
a_{s}\frac{d^{N-1}}{dx^{N-1}}x^{s},...,\left(  2k\right)  !%
{\displaystyle\sum_{s=N}^{\infty}}
a_{s}\frac{d^{N-2k}}{dx^{N-2k}}x^{s}\right)
\]

Using continuity and linearity arguments we see that%
\[
A_{N,2k}\left(  x\right)  =%
{\displaystyle\sum_{s=N}^{\infty}}
a_{s}A_{N,2k,s}\left(  x\right)  ,
\]
where the vectors $A_{N,2k,s}\left(  x\right)  $ are defined by
\[
A_{N,2k,s}\left(  x\right)  :=\left(  \frac{d^{N}}{dx^{N}}x^{s},1!\frac
{d^{N-1}}{dx^{N-1}}x^{s},...,\left(  2k\right)  !\frac{d^{N-2k}}{dx^{N-2k}%
}x^{s}\right)
\]

In order to prove for every $x\in\mathbb{R}_{+}$ inequality
(\ref{POSITIVE-whole-AXIS}), we have to show that for any $p\in\mathbb{R}%
^{k+1}$ with $p\neq0$ we have the inequality:%
\[
Q_{1}\left(  p\right)  :=%
{\displaystyle\sum_{i,j=0}^{k}}
p_{i}p_{j}b_{i+j}\left(  x\right)  >0.
\]
Clearly we have
\[
Q_{1}\left(  p\right)  =%
{\displaystyle\sum_{s=N}^{\infty}}
a_{s}\cdot%
{\displaystyle\sum_{i,j=0}^{k}}
p_{i}p_{j}A_{N,2k,s;i+j}\left(  x\right)  ,
\]
where $A_{N,2k,s;m}\left(  x\right)  $ is the $m$th entry of the vector
$A_{N,2k,s}\left(  x\right)  .$

For $s\geq N$ we have%
\[
\frac{d^{N}}{dx^{N}}x^{s}=\frac{s!}{\left(  s-N\right)  !}x^{s-N}%
\]
hence
\[
A_{N,2k,s}\left(  x\right)  =x^{s-N}\left(  \frac{s!}{\left(  s-N\right)
!},\frac{s!}{\left(  s-N+1\right)  !}x,...,\left(  2k\right)  !\frac
{s!}{\left(  s-N+2k\right)  !}x^{2k}\right)  .
\]

We extract the factor $x^{\left(  s-N\right)  \left(  k+1\right)  }$ and
obtain the vectors $\widetilde{A}_{N,k,s}\left(  x\right)  ,$ namely,%
\[
A_{N,2k,s}\left(  x\right)  =x^{s-N}\widetilde{A}_{N,2k,s}\left(  x\right)  .
\]
where
\[
\widetilde{A}_{N,2k,s}\left(  x\right)  :=\left(  \frac{s!}{\left(
s-N\right)  !},\frac{s!}{\left(  s-N+1\right)  !}x,...,\left(  2k\right)
!\frac{s!}{\left(  s-N+2k\right)  !}x^{2k}\right)
\]

Obviously, due to $x>0,$ we have $x^{s-N}>0.$ Since $a_{s}\geq0,$ if we prove
that for every $x\geq0$ the vectors $A_{N,k,s}\left(  x\right)  $ satisfy the
analogs to inequalities (\ref{POSITIVE-HALF-AXIS}) and
(\ref{POSITIVE-whole-AXIS}), then our statement will be proved. Hence, we have
to show\textbf{ }that the vectors\textbf{ }$\widetilde{A}_{N,2k,s}\left(
x\right)  $ satisfy analogs to inequalities (\ref{POSITIVE-HALF-AXIS}) and
(\ref{POSITIVE-whole-AXIS}).

By a classical result (see \cite{Horn-Johnson}), it is sufficient to prove
that the \textbf{Hankel matrix} (we apply the usual indexing of Hankel
matrices) defined by \textbf{ }
\[
\widetilde{H}_{N,k,s}\left(  x\right)  =\left(
\begin{array}
[c]{cccc}%
\frac{s!}{\left(  s-N\right)  !} & \frac{s!}{\left(  s-N+1\right)  !}x &
.... & k!\frac{s!}{\left(  s-N+k\right)  !}x^{k}\\
\frac{s!}{\left(  s-N+1\right)  !}x & ... & k!\frac{s!}{\left(  s-N+k\right)
!}x^{k} & \left(  k+1\right)  !\frac{s!}{\left(  s-N+k+1\right)  !}x^{k+1}\\
& ... &  & \\
k!\frac{s!}{\left(  s-N+k\right)  !}x^{k} &  &  & \left(  2k\right)
!\frac{s!}{\left(  s-N+2k\right)  !}x^{2k}%
\end{array}
\right)
\]
is \textbf{strictly} positive-definite for $s>N$; it is easy to see that for
$s=N$ we have only positive semi-definite. For proving the strict
positive-definiteness, we apply the \textbf{determinant criterion} (see
\cite{Horn-Johnson}), i.e. we have to see that for every fixed $s>N$
\[
\det\widetilde{H}_{N,k,s}\left(  x\right)  >0\qquad\text{ for all }k\geq1.
\]
By extracting the powers of $x$ first by columns and then by rows, we obtain
\begin{align*}
&  \det\widetilde{H}_{N,k,s}\left(  x\right) \\
&  =x^{k\left(  k+1\right)  }\det\left(
\begin{array}
[c]{cccc}%
\frac{s!}{\left(  s-N\right)  !} & \frac{s!}{\left(  s-N+1\right)  !} & .... &
k!\frac{s!}{\left(  s-N+k\right)  !}\\
\frac{s!}{\left(  s-N+1\right)  !} & ... & k!\frac{s!}{\left(  s-N+k\right)
!} & \left(  k+1\right)  !\frac{s!}{\left(  s-N+k+1\right)  !}\\
& ... &  & \\
k!\frac{s!}{\left(  s-N+k\right)  !} &  &  & \left(  2k\right)  !\frac
{s!}{\left(  s-N+2k\right)  !}%
\end{array}
\right)  ,
\end{align*}
hence, we see that $x^{k\left(  k+1\right)  }>0$ for every $x$ since $k\left(
k+1\right)  $ is an even number. Further, by putting
\[
S:=s-N
\]
we obtain the representation%
\begin{align*}
&  \det\widetilde{H}_{N,k,s}\left(  x\right)  =\\
&  =\frac{x^{k\left(  k+1\right)  }\left(  s!\right)  ^{k+1}}{\left(
S!\right)  ^{k}}\det\left(
\begin{array}
[c]{cccc}%
\frac{S!}{S!} & \frac{S!}{\left(  S+1\right)  !} & .... & k!\frac{S!}{\left(
S+k\right)  !}\\
\frac{S!}{\left(  S+1\right)  !} & ... & k!\frac{S!}{\left(  S+k\right)  !} &
\left(  k+1\right)  !\frac{S!}{\left(  S+k+1\right)  !}\\
& ... &  & \\
k!\frac{S!}{\left(  S+k\right)  !} &  &  & \left(  2k\right)  !\frac
{S!}{\left(  S+2k\right)  !}%
\end{array}
\right)  .
\end{align*}
The positivity of $\det\widetilde{H}_{N,k,s}\left(  x\right)  $ follows now
directly from the theorem of Chammam\textbf{ }\cite{Chammam2019},
\cite{Chu2021}, see below the Remark (\ref{RemarkChammam}), where we have put
\[
\alpha=1,\qquad1+\beta+\alpha=S+1\Longrightarrow\beta=S-1
\]
Note that by our assumptions, we have $S\geq1$. Indeed, we have
\[
\left(  1\right)  _{n}=n!\qquad\text{for }n\geq0.
\]

2. Further, we have to prove inequality (\ref{POSITIVE-HALF-AXIS}) for every
$x\in\mathbb{R}_{+},$ i.e. we have to consider the quadratic forms
\[
Q_{2}:=%
{\displaystyle\sum_{i,j=0}^{k}}
p_{i}p_{j}\widehat{c}_{i+j+1}\left(  x\right)  .
\]
Similar to the previous case, we may extract the powers of $x,$ first the
powers $x^{s-N+1}$ from the vector $A_{N,2k,s}^{\prime}\left(  x\right)  $
(which is equal to the vector $A_{N,2k,s}\left(  x\right)  $ with the first
element dropped) and then in the corresponding Hankel determinant $\det\left(
\widetilde{H}_{N,k,s}^{\prime}\left(  x\right)  \right)  $ we extract the
positive factor $x^{k\left(  k+1\right)  }$ . Since we have
\begin{align*}
&  \det\left(
\begin{array}
[c]{cccc}%
\frac{1}{\left(  S+1\right)  !} & \frac{2}{\left(  S+2\right)  !} & .... &
k!\frac{1}{\left(  S+k\right)  !}\\
... &  & k!\frac{1}{\left(  S+k\right)  !} & \left(  k+1\right)  !\frac
{1}{\left(  S+k+1\right)  !}\\
... &  &  & \\
k!\frac{1}{\left(  S+k\right)  !} &  &  & \left(  2k-1\right)  !\frac
{1}{\left(  S+2k-1\right)  !}%
\end{array}
\right)  \\
&  =\frac{1}{\left(  \left(  S+1\right)  !\right)  ^{k}}\det\left(
\begin{array}
[c]{ccccc}%
\frac{\left(  S+1\right)  !}{\left(  S+1\right)  !} & \frac{\left(
S+1\right)  !}{\left(  S+2\right)  !}2 & \cdot\cdot\cdot &  & k!\frac{\left(
S+1\right)  !}{\left(  S+k\right)  !}\\
... &  &  & k!\frac{\left(  S+1\right)  !}{\left(  S+k\right)  !} & \left(
k+1\right)  !\frac{\left(  S+1\right)  !}{\left(  S+k+1\right)  !}\\
... &  &  &  & \\
k!\frac{\left(  S+1\right)  !}{\left(  S+k\right)  !} &  &  &  & \left(
2k-1\right)  !\frac{\left(  S+1\right)  !}{\left(  S+2k-1\right)  !}%
\end{array}
\right)
\end{align*}
we obtain the determinants
\begin{align*}
&  \det\left(  \widetilde{H}_{N,k,s}^{\prime}\left(  x\right)  \right)  \\
&  =\frac{x^{k\left(  k+1\right)  }\left(  s!\right)  ^{k+1}}{\left(
S!\right)  ^{k}}\times\det\left(
\begin{array}
[c]{cccc}%
\frac{1}{\left(  S+1\right)  !} & \frac{2}{\left(  S+2\right)  !} & .... &
k!\frac{1}{\left(  S+k\right)  !}\\
... &  & k!\frac{1}{\left(  S+k\right)  !} & \left(  k+1\right)  !\frac
{1}{\left(  S+k+1\right)  !}\\
... &  &  & \\
k!\frac{1}{\left(  S+k\right)  !} &  &  & \left(  2k-1\right)  !\frac
{1}{\left(  S+2k-1\right)  !}%
\end{array}
\right)  \\
&  =\frac{x^{k\left(  k+1\right)  }\left(  s!\right)  ^{k+1}}{\left(
S!\right)  ^{k}\left(  \left(  S+1\right)  !\right)  ^{k}}\det\left(
\begin{array}
[c]{ccccc}%
\frac{\left(  S+1\right)  !}{\left(  S+1\right)  !} & \frac{\left(
S+1\right)  !}{\left(  S+2\right)  !}2 & \cdot\cdot\cdot &  & k!\frac{\left(
S+1\right)  !}{\left(  S+k\right)  !}\\
... &  &  & k!\frac{\left(  S+1\right)  !}{\left(  S+k\right)  !} & \left(
k+1\right)  !\frac{\left(  S+1\right)  !}{\left(  S+k+1\right)  !}\\
... &  &  &  & \\
k!\frac{\left(  S+1\right)  !}{\left(  S+k\right)  !} &  &  &  & \left(
2k-1\right)  !\frac{\left(  S+1\right)  !}{\left(  S+2k-1\right)  !}%
\end{array}
\right)
\end{align*}
and in every term above we obtain the following expressions:\
\[
\frac{\left(  j+1\right)  !}{\left(  S+2\right)  \cdot\cdot\cdot\left(
S+j+1\right)  }=\frac{\left(  2\right)  _{j}}{\left(  S+2\right)  _{j}}%
\qquad\text{for }j\geq1
\]
Hence, the determinant above is obtained in an explicit form again from the
result of \cite{Chammam2019}, \cite{Chu2021}, see Remark \ref{RemarkChammam}
below, by putting
\[
\alpha=1,\qquad1+\alpha+\beta=S.
\]

3. Somewhat more complicated is the case of the compact interval $\left[
a,b\right]  ,$ where we have taken for simplicity the interval $\left[
0,1\right]  .$ Since $x\in\left(  0,1\right)  $ implies $x>0,$ we prove
inequality (\ref{POSITIVE-HALF-AXIS}) as in the previous section.

Further, for every $x\in\left(  0,1\right)  $, we have to prove inequality
(\ref{POSITIVE-INTERVAL}), hence we consider the third quadratic form
\[
Q_{3}:=%
{\displaystyle\sum_{i,j=0}^{k}}
p_{i}p_{j}\left(  b_{i+j}\left(  x\right)  -b_{i+j+1}\left(  x\right)
\right)
\]
As in the previous paragraphs, it is not difficult to see that for $Q_{3}$ we
have the elements in the vector $A_{N,k,s}\left(  x\right)  $ equal to\
\begin{align*}
c_{k}\left(  x\right)   &  :=x^{k}\frac{k!}{\left(  S+k\right)  !}%
-x^{k+1}\frac{\left(  k+1\right)  !}{\left(  S+k+1\right)  !}=\frac{x^{k}%
k!}{\left(  S+k\right)  !}\left(  1-x\frac{k+1}{S+k+1}\right) \\
&  =\frac{k!}{\left(  S+k\right)  !}\frac{S+k+1-x\left(  k+1\right)  }%
{S+k+1}=\frac{k!}{\left(  S+k+1\right)  !}\left(  S+\left(  1-x\right)
\left(  k+1\right)  \right)
\end{align*}
Now, we may prove easily that
\[
Q_{3}\left(  x\right)  :=%
{\displaystyle\sum_{0\leq i,j\leq k}}
p_{i}p_{j}c_{i+j}\left(  x\right)  \geq0
\]
Indeed, we have
\begin{align*}
Q_{3}\left(  x\right)   &  =S%
{\displaystyle\sum_{0\leq i,j\leq k}}
p_{i}p_{j}\frac{\left(  i+j\right)  !}{\left(  S+i+j+1\right)  !}\\
&  +\left(  1-x\right)
{\displaystyle\sum_{0\leq i,j\leq k}}
p_{i}p_{j}\frac{\left(  i+j+1\right)  !}{\left(  S+i+j+1\right)  !}%
\end{align*}
and we may apply to both sums the formula of Chammam\textbf{ }%
\cite{Chammam2019}, \cite{Chu2021}, see Remark \ref{RemarkChammam}, to prove
that the determinants of the corresponding Hankel matrix are $\geq0,$ which is
equivalent to the positive-definiteness of the matrices. Note that we use
essentially the fact that $0\leq x\leq1.$

This ends the proof.
\end{proof}

\begin{proposition}
The inequality $\Lambda\geq0,$ i.e. $\lambda_{j}\geq0$ for all $j,$ implies
\[
\Phi_{N}^{\left(  j\right)  }\left(  x\right)  \geq0\text{ for all }x>0\text{
and all }j\geq N.
\]

\end{proposition}

Indeed, this is proved by the Taylor expansion of the function $\Phi_{N}$, in
paper \cite{NATO-Paper}, see also
\cite{AldazKounchevRender-ConstructiveApproximation-2008}.

\begin{remark}
\label{RemarkChammam} For the proofs above we used the results of Chammam
\cite{Chammam2019} about determinants of Hankel matrices, see also Chu
\cite{Chu2021}.. For numbers $x\in\mathbb{C}$ and a natural number $m\geq1,$
we have the notations (as in \cite{Chu2021}):\
\[
\left(  x\right)  _{0}:=1,\qquad\left(  x\right)  _{m}:=x\left(  x+1\right)
\cdot\cdot\cdot\left(  x+m-1\right)  \qquad\text{for }m\geq1
\]
hence, obviously,
\[
\left(  1\right)  _{m}=m!
\]

Then, according to formula $1.1$ in \cite{Chu2021}, we have the main result of
Chammam:\
\[
\det_{0\leq i,j\leq m}\left[  \frac{\left(  \alpha\right)  _{i+j}}{\left(
1+\alpha+\beta\right)  _{i+j}}\right]  =\prod_{k=0}^{m}\frac{k!\left(
\alpha\right)  _{k}\left(  1+\beta\right)  _{k}}{\left(  \alpha+\beta
+k\right)  _{k}\left(  1+\alpha+\beta\right)  _{2k}}.
\]

\end{remark}

\begin{remark}
Obviously, the results in Theorem \ref{THM-Main} may be generalized to
arbitrary compact interval $\left[  a,b\right]  $ which is a subset of
$\left[  0,\infty\right)  $ and also for non-compact $\left[  a,b\right]
\subset\left[  0,\infty\right)  $ with $b=\infty.$
\end{remark}

Now we apply the above Theorem to the Exponential Moment problem.

\begin{theorem}
\label{THM-Moments} Assume that the non-negative measure $d\mu$ is given on
the interval (compact or not) $\left[  a,b\right]  \subset\left[
0,\infty\right)  .$ We put
\[
\widehat{c}_{j}:=%
{\displaystyle\int_{a}^{b}}
b_{j}\left(  x\right)  d\mu\left(  x\right)  \qquad\text{for }j=0,1,...,N
\]
where $b_{j}\left(  x\right)  $ are given by (\ref{bj-DEFINITION}).

Then there exists a non-negative measure $d\nu$ defined on the interval
$\left[  a,b\right]  $ which solves the classical Moment problem
\[%
{\displaystyle\int_{a}^{b}}
t^{j}d\nu\left(  t\right)  =\widehat{c}_{j}\qquad\text{for }j=0,1,...,N.
\]

\end{theorem}

The proof follows from Theorem \ref{THM-Main}, by the Theorems of M. Riesz and
Stieltjes provided in \cite{Akhiezer}, see above Remark \ref{REMARK-1}.


\begin{thebibliography}{99}                                                                                               %


\bibitem {Akhiezer}Akhiezer, N.I.: \emph{The problem of moments and some
related questions in analysis,} Oliver \& Boyd, Edinburgh, 1965.(Transl. from
Russian ed. Moscow 1961).

\bibitem {AldazKounchevRender-ConstructiveApproximation-2008}Aldaz, J. M.;
Kounchev, O.; Render, H.: \emph{Bernstein operators for exponential
polynomials.} Constr. Approx. 29 (2009), no. 3, 345--367.

\bibitem {Berenstein}Berenstein C. A.; Gay, R.: \emph{Complex analysis and
special topics in harmonic analysis}, Springer-Verlag, New York, 1995.

\bibitem {Chu2021}Chu, W., \emph{Hankel determinants of factorial fractions},
Bull. Aust. Math. Soc. Vol. 105 , Issue 1, 2022, pp. 46 - 57

\bibitem {Chammam2019}Chammam, W., \emph{Generalized harmonic numbers, Jacobi
numbers and a Hankel determinant evaluation,} Integral Transforms Spec. Funct.
30(7) (2019), 581--593.

\bibitem {Gustafsson-Putinar}Gustafsson, B. and M. Putinar, \emph{Hyponormal
Quantization of Planar Domains. Exponential Transform in Dimension Two},
Springer, $2017.$

\bibitem {Horn-Johnson}Horn, Roger A.; Johnson, Charles R. (2013).
\emph{Matrix Analysis} (2nd ed.). Cambridge University Press.

\bibitem {kounchevBOOK}Kounchev, O.I.: \emph{Multivariate Polysplines.
Applications to Numerical and Wavelet Analysis}, Academic Press, London--San
Diego, 2001.

\bibitem {NATO-Paper}Kounchev, O. and H. Render, \emph{On a new method for
controlling exponential processes, }Proc. Conference\emph{ }"Scientific
Support for the Decision Making in the Security Sector" Edited by O. Kounchev,
et al. , IOS Press, 2007; online at \textbf{ }https://arxiv.org/abs/0901.0327.

\bibitem {KoReArkiv2010}Kounchev, O.; Render, H.: \emph{A moment problem for
pseudo-positive definite functionals.} Ark Mat 48, 97--120 (2010).

\bibitem {KoRe13}Kounchev, O.; Render, H.: \emph{Polyharmonic functions of
infinite order on annular regions,} T\^{o}hoku Math. J. 65 (2013), 199--229.

\bibitem {KoReInter2019}Kounchev, O.; Render, H.: \emph{Interpolation of data
functions on parallel hyperplanes,} J. Approx. Theory 246 (2019), 43--61.

\bibitem {KounchevRender2021JCAM}Kounchev, O.; Render, H.: \emph{Error
estimates for interpolation with piecewise exponential splines of order two
and four}, J. Comput. App. Math., Vol. 391, 1 August $2021$, $113464.$

\bibitem {KounchevRenderTsachevBIT}Kounchev, O.; Render, H.; Tsachev, Ts.:
\emph{On a class of }$L-$\emph{splines of order }$4$\emph{: fast algorithms
for interpolation and smoothing}, BIT Numerical Mathematics, volume 60, pages
879--899 (2020).

\bibitem {McCa91}McCartin, B.J.: \emph{Theory of exponential splines,} J.
Approx. Theory 66 (1991), 1--23.

\bibitem {Markov-1886}Markov, A., \emph{Sur une question de maximum et de
minimum proposee par M. Tchebycheff}, Acta Math. 9 (1986/1987), 57-70.

\bibitem {KreinNudelman}Krein, M. and A. Nudelman, \emph{The Markov Moment
Problem and Extremal Problems}, Translation of AMS, 1977.

\bibitem {RamsaySilverman}Ramsay, J.O.; Silverman, B. W.: \emph{Functional
Data Analysis}, Springer Verlag, Second Edition 2005.

\bibitem {Schumaker}Schumaker, L.L.: \emph{Spline Functions: Basic Theory},
Interscience, New York, 1981.

\bibitem {Unser}Unser, M.; Blu, T.: \emph{Cardinal Exponential Splines: Part I
-- Theory and Filtering Algorithms,} IEEE Transactions on Signal Processing,
53 (2005), 1425--1438.
\end{thebibliography}
\end{document}